\documentclass[a4paper,10pt]{article}

\usepackage[english]{babel}
\usepackage{german, umlaute, amsmath, amssymb , latexsym, theorem,
  epic, epsfig, rotating}

\newcommand{\an}{\ensuremath{a_{n}}}

 \newcommand{\qn}{\ensuremath{q_{n}}}

\newcommand{\nnull}{\ensuremath{n_{0}}}

\newcommand{\neins}{\ensuremath{n_{1}}}

\newcommand{\atel}{\ensuremath{\frac{1}{a}}}

\newcommand{\halb}{\ensuremath{\frac{1}{2}}}

\newcommand{\alphlist}{\begin{list}{(\alph{enumi})}{\usecounter{enumi}}}
\newcommand{\romanlist}{\begin{list}{(\roman{enumi})}{\usecounter{enumi}}}
\newcommand{\listend}{\end{list}}

\newcommand{\ld}{\ensuremath{,\ldots,}}

\newcommand{\ra}{\ensuremath{\rightarrow}}

\newcommand{\follows}{\ensuremath{\Rightarrow}}

\newcommand{\N}{\ensuremath{\mathbb{N}}} 
\newcommand{\R}{\ensuremath{\mathbb{R}}}

\newcommand{\Q}{\ensuremath{\mathbb{Q}}}

\newcommand{\kreis}{\ensuremath{\mathbb{T}^{1}}}

\newcommand{\nfolge}[1]{\ensuremath{(#1)_{n\in\mathbb{N}}}}

\newcommand{\qed}{{\raggedleft $\Box$ \\}}
\newcommand{\proof}{\textit{Proof:}}

\newcommand{\Tin}{\ensuremath{T^{-1}}}

\newcommand{\incap}{\ensuremath{\bigcap_{i=1}^n}}

\newcommand{\thx}{\ensuremath{(\theta,x)}}
\newcommand{\thom}{\ensuremath{\theta + \omega}}

\newcommand{\clPhiplus}{\ensuremath{\overline{\Phi^+}}}

\newcommand{\phith}{\ensuremath{\varphi(\theta)}}

\newcommand{\graphpt}{\ensuremath{(\theta,\varphi(\theta))}}

\newcommand{\Tth}{\ensuremath{T_{\theta}}}

\newcommand{\Tthx}{\ensuremath{T_{\theta}(x)}}

\theoremstyle{break}
\newtheorem{definition}{Definition}[section]
\newtheorem{thm}[definition]{Theorem}

\newtheorem{lem}[definition]{Lemma}
\newtheorem{cor}[definition]{Corollary}

\theorembodyfont{\rmfamily}
\newtheorem{bem}[definition]{Remark}    
\newtheorem{example}[definition]{Example}

\numberwithin{equation}{section}

\title{On the structure of strange nonchaotic attractors in pinched
  skew products} \author{Tobias H. J\"ager \\ 
  \small Mathematisches
  Institut,  \\
  \small Friedrich-Alexander-Universit\"at Erlangen-N\"urnberg,
  Germany \\ \small MSC Classification Numbers: 37C70, 37C60, 37H15 \\
}

\begin{document}
\setlength{\topmargin}{-1cm}
\setlength{\oddsidemargin}{-0.03\textwidth}  

\maketitle

\begin{abstract}
  The existence of non-continuous invariant graphs (or strange
  non-chaotic attractors) in quasiperiodically forced systems has
  generated great interest, but there are still very few rigorous
  results about the properties of these objects. In particular it is
  not known whether the topological closure of such graphs is
  typically a filled-in set, i.e.\ consist of a single interval on
  every fibre, or not. We give a positive answer to this question for
  the class of so-called pinched skew products, where non-continuous
  invariant graphs occur generically, provided the rotation number on
  the base is diophantine and the system satisfies some additional
  conditions.  For typical parameter families these conditions
  translate to a lower bound on the parameter. On the other hand, we
  also construct examples where the non-continuous invariant
  graphs contain isolated points, such that their topological closure
  cannot be filled in.
\end{abstract}


\section{Introduction}

Generally spoken, for a quasiperiodically forced monotone interval map
invariant graphs play the same role as fixed points for unperturbed
maps. However, in contrast to fixed points there are two different
types of invariant graphs. In the simpler case, such a graph is
continuous. If in addition it has a negative Lyapunov exponent, a lot
of conclusions about its regularity, structural stability etc.\ can be
drawn (see \cite{stark:1999,sturman/stark:2000}). More difficult, but
also more interesting, is the case of non-continuous invariant graphs.
When their Lyapunov exponent is negative, such graphs are often
refered to as strange non-chaotic attractors or SNA. Their existence
has recently received great attention in theoretical physics, and
consequently there are a lot of numerical studies about the topic
(\cite{prasad/negi/ramaswamy:2001} gives a good overview and further
reference). On the other hand, rigorous results are still few.
The existence of SNA has only been proved for rather specific systems
(see \cite{herman:1983} for quasiperiodically driven Moebius
transformations, \cite{keller:1996,bezhaeva/oseledets:1996} for
pinched skew products), and even in these cases there are a lot of
open questions regarding their further properties. It is one of these
questions we want to address here.

\ \\
\textbf{Quasiperiodically forced monotone interval maps.} In order to
give a precise formulation of the problem we want to study, we
need a few basic definitions. First of all, a
\textit{quasiperiodically forced monotone interval map} is a
continuous map of the form
\begin{equation}
  \label{eq:qpfmim}
  T : \kreis \times I \ra \kreis \times I \ , \ \thx \mapsto
  (\thom,\Tthx) 
\end{equation}
where $I \subset \R$ is a compact interval and all the \textit{fibre
  maps} $\Tth$ are monotonically increasing on $I$. An
\textit{invariant graph} for such a system $T$ is a function $\varphi
: \kreis \ra I$ which satisfies
\begin{equation}
  \label{eq:invgraph}
  \Tth(\phith) \ = \ \varphi(\thom) \ \ \forall \theta \in \kreis .
\end{equation}
This equation implies that the corresponding point set $\Phi := \{
\graphpt \mid \theta \in \kreis \}$ is forward invariant, i.e.\ 
$T(\Phi) = \Phi$ (not necessarily $\Tin(\Phi) = \Phi$, as we
did not assume strict monotonicity).
On the other hand, whenever $K$ is a compact and forward invariant
set we can define an invariant graph by
\begin{equation}
  \label{eq:boundinggraph}
  \varphi^+(\theta) \ := \ \max\{ x \in I \mid \thx \in K \}
\end{equation}
(the invariance being a direct
consequence of the monotonicity of the fibre maps).
Further, as $K$ is compact $\varphi^+$ will be upper semi-continuous.
In the same way a lower semi-continuous invariant graph $\varphi^-$
can be defined via the minimum.  Particularly interesting for our
purposes will be the case where $K$ is the \textit{global attractor}
of the system, defined as $K := \bigcap_{n=0}^\infty T^n(\kreis \times
I)$. The corresponding graphs $\varphi^+$ and $\varphi^-$ will then be
called the \textit{upper and lower bounding graphs} of the system $T$.
(All of this is throughly discussed in \cite{glendinning:2002},
\cite{stark:2003} or \cite{jaeger:2003}.)

\ \\
\textbf{The question.} Suppose $\varphi$ is an invariant graph
which is not continuous.  Its topological closure $\overline{\Phi}$ is
a compact and invariant set, again bounded from above and below by
invariant graphs $\varphi^+$ and $\varphi^-$. If we let
$[\varphi^-,\varphi^+] := \{ \thx \mid \varphi^-(\theta) \leq x \leq
\varphi^+(\theta) \}$, then this is a compact invariant set as well and 
surely $\overline{\Phi^+} \subseteq [\varphi^-,\varphi^+]$.
But is $\clPhiplus = [\varphi^-,\varphi^+]$?

This question was already asked by M.\ Herman in \cite{herman:1983}
(Section 4.14) for certain quasiperiodically forced
Moebius-transfor\-mations, and then repeated a number of times in
different situations (see
\cite{keller:1996,glendinning:2002,stark:2003}).

\ \\ 
\textbf{Results.} Here, we address the
problem in the setting of so-called pinched skew products first
introduced in \cite{grebogi/ott/pelican/yorke:1984}. Their particular
structure allows to prove the existence of non-continuous invariant
graphs by a few simple and elegant arguments (see
\cite{keller:1996,bezhaeva/oseledets:1996}, a slight generalization
can be found in \cite{glendinning:2002}). Often, such systems are
given by (\ref{eq:qpfmim}) with $I = [0,L]$, $L>0$ and fibre maps
\begin{equation} \label{eq:pinchedskews}
    \Tthx = \alpha g(\theta)f(x) \ .
\end{equation}
where $f : \R^+ \ra \R^+$ is monotonically increasing with $f(0)=0$,
$g : \kreis \ra \R^+$ has exactly one zero and $\alpha$ is a positive
parameter. A typical example would be $f(x)= \tanh(x)$ and $g(\theta)
= |\sin(\pi\theta)|$. Note that $f(0)=0$ implies that the lower
bounding graph is always $\varphi^- \equiv 0$. In short terms, our
main result can now be stated as follows (see Cor.\ 
\ref{cor:fgsystem}):
\begin{quote}
  Suppose $f$ and $g$ satisfy some mild conditions concerning their
  geometry and regularity (specified before Cor.\ \ref{cor:fgsystem})
  and the rotation number $\omega$ is of diophantine type. Then for
  all sufficiently large parameters $\alpha$ the upper bounding graph
  $\varphi^+$ of the system $T$ defined by (\ref{eq:pinchedskews})
  will have the property
$\clPhiplus = [0,\varphi^+]$.
\end{quote}
As mentioned, we also give examples where the topological closure of
the upper bounding graph cannot be filled in because the graph
contains isolated points. For this we use a certain growth condition
on the coefficients $a_n$ of the continued fraction expansion of the
rotation number $\omega$ on the base, namely $\sum_{i=1}^\infty
\frac{1}{a_i} < \infty$. However, this still allows $\omega$ to be
diophantine, showing that the additional assumptions on the system we
use to derive the above result are indeed crucial, and cannot be
neglected.

\ \\
\textbf{Overview.} Section \ref{Pinchedskews} briefly sketches the
arguments used to establish the existence of non-continuous invariant
graphs in pinched skew products. In order to gain more insight about
their creation and structure, Section \ref{Boundarylines} introduces
the \textit{iterated upper boundary lines}. This sequence of
continuous graphs converges monotonically decreasing to the upper
bounding graph, and their shape can be controlled to some extent by
putting additional restrictions on the geometry of the system.  In
Section \ref{Mainresult} this is then used to derive the mentioned
result, and Section \ref{Counterexamples} contains the construction of 
the counterexamples.

\ \\
\textbf{Acknowlegements.} I thank Gerhard Keller for introducing me to
the problem, and for many helpful remarks and discussions concerning
this manuscript. 


\section{Pinched skew products} \label{Pinchedskews}

Suppose a map $T$ as in (\ref{eq:qpfmim}) with $I = [0,L]$ for some
$L>0$ additionally satisfies
\begin{eqnarray}
  & \bullet & \Tth(0) = 0 \ \forall \theta\in\kreis \hspace{10,4eM}
  \textrm{(invariant 0-line)} \hspace{7eM} \label{eq:0line} \\
  & \bullet & T_{\theta} \equiv 0 \ \textrm{ for some } \ \theta
  \in \kreis \hspace{7,5eM} \textrm{(pinching)}
  \label{eq:pinching}
\end{eqnarray}
As mentioned, (\ref{eq:0line}) ensures that the lower bounding graph
of such a system is always $\varphi^- \equiv 0$. Further
(\ref{eq:pinching}) implies that any other invariant graph must equal
0 on a dense subset of $\kreis$, namely the forward orbit of
$\theta$. Thus, apart from $\varphi^-$ no other invariant graph can
be continuous. 

Now suppose all the fibre maps $\Tth$ are differentiable and denote
the derivative by $D\Tth$. Then by
\begin{equation}
  \label{eq:lyapunov}
  \lambda(\varphi) \ := \ \int_{\kreis} \log D\Tth(\varphi(\theta)) \
  d\theta 
\end{equation}
the Lyapunov exponent of an invariant graph $\varphi$ can be defined.
It is easy to show that whenever $T$ satisfies some mild conditions
(namely when $\theta \mapsto \inf_{x \in I}D\Tthx$ is integrable), the
Lyapunov exponent of the upper bounding graph cannot be strictly
positive (e.g.\ Lemma 3.5 in \cite{jaeger:2003}). For
systems of the form (\ref{eq:pinchedskews}), a simple computation
yields
\begin{equation}  \label{eq:snacondition}
    \lambda(0) \ = \ \log f'(0) + \log \alpha + \int_{\kreis} \log
    g(\theta) \ d\theta
\end{equation}
whenever $f'(0) > 0$ and $\log g$ is integrable. This is surely
positive for sufficiently large $\alpha$. Thus the upper bounding
graph $\varphi^+$ cannot be equal to $\varphi^-$ and must therefore be 
non-continuous in such systems.

Even if they are rather degenerate in some sense, the fact that the
existence of strange non-chaotic attractors can be established so
easily makes pinched skew products an ideal setting for studying their
further properties. The results obtained here may then at least give
hints about more general systems.  This is further supported by the
fact that although invertible systems lack the pinched structure,
their minimal sets will still be pinched sets (i.e.\ consist of only
one point on a dense set of fibres, see \cite{stark:2003}), and
therefore have some similarities with the global attractors of pinched
skew products.

\begin{figure}[h!]
\noindent
\begin{minipage}[t]{\linewidth}
  \epsfig{file=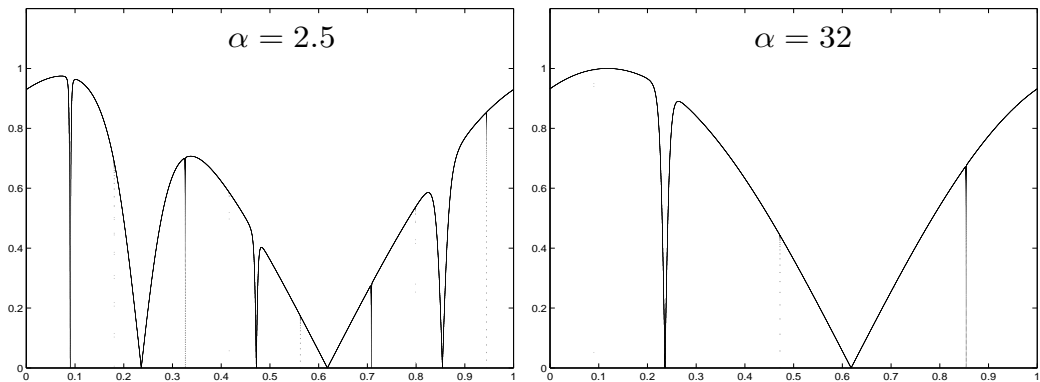, clip=, width=\linewidth,
    height=0.25\textheight}
         \caption{\small Semi-continuous
          upper bounding graphs in the parameter family $\thx \mapsto
          (\thom, \tanh(\alpha x) \cdot |\sin(\pi\theta)|)$. Here
          $\omega$ is the golden mean and $\alpha$ is 5 and 32,
          respectively. 
         }             \label{fig:sna} 
\end{minipage}
\end{figure}


\section{The shape of the iterated upper boundary lines}
\label{Boundarylines}

The argument sketched in the last section gives the non-continuity of
the upper bounding graph and ensures that its topological closure
contains the 0-line, but apart from that it offers little further
information. To explain the shape of such SNA as depicted in Figure
\ref{fig:sna}, recall that the upper bounding graph was defined via
the global attractor $K = \bigcap_{n \in \N} T^n(\kreis \times I)$.
The sets $K_n := \incap T^i(\kreis \times I)$ are bounded above by the
iterates of the upper boundary line $\kreis \times \{ L \}$. To be
more precise, define a monotonically decreasing sequence of continuous
graphs by $\varphi_0 :\equiv L$ and
$\varphi_{n+1}(\theta) := T_{\theta-\omega}(\varphi_n(\theta-\omega))$
(alternatively $\varphi_n(\theta) := T^n_{\theta-n\omega}(L)$). Then
$K_n = [0,\varphi_n]$, and the graphs $\varphi_n$ converge pointwise
and monotonically decreasing to $\varphi^+$. It therefore seems reasonable
to hope that by understanding the behaviour of the $\varphi_n$ we can
gain more insight about their limit object.

\begin{figure}[h!]
\noindent
\begin{minipage}[t]{\linewidth}
  \epsfig{file=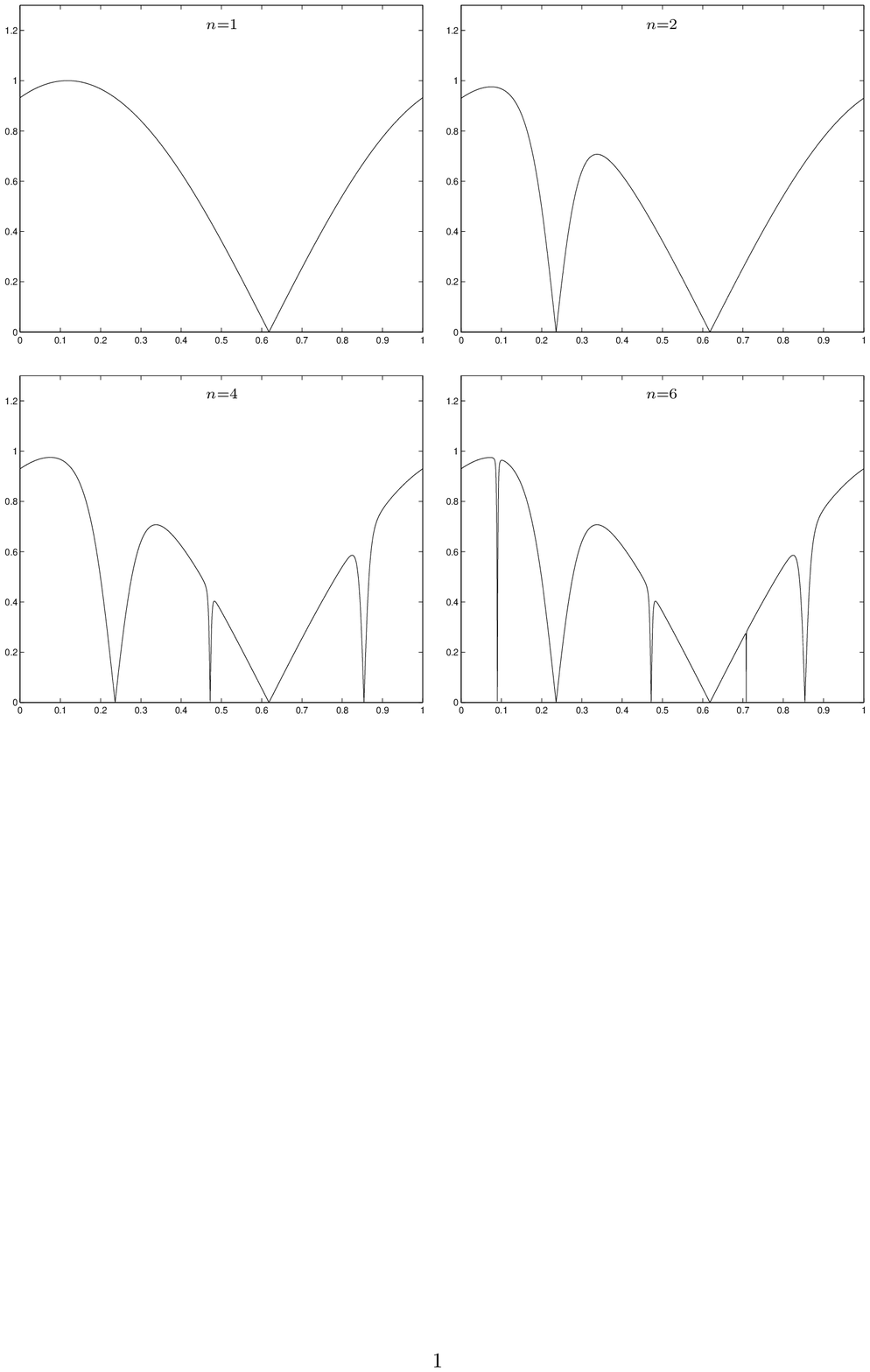, clip=, width=\linewidth,
    height=0.5\textheight} \caption{\small The iterated upper boundary
    lines $\varphi_n$ of the system $\thx \mapsto (\thom,\tanh(5x)
    \cdot |\sin(\pi\theta)|)$. Again $\omega$ is the golden mean.  }
  \label{fig:bounding}
\end{minipage}
\end{figure}

\begin{figure}[h!]   
\noindent
\begin{minipage}[t]{\linewidth} 
 \epsfig{file=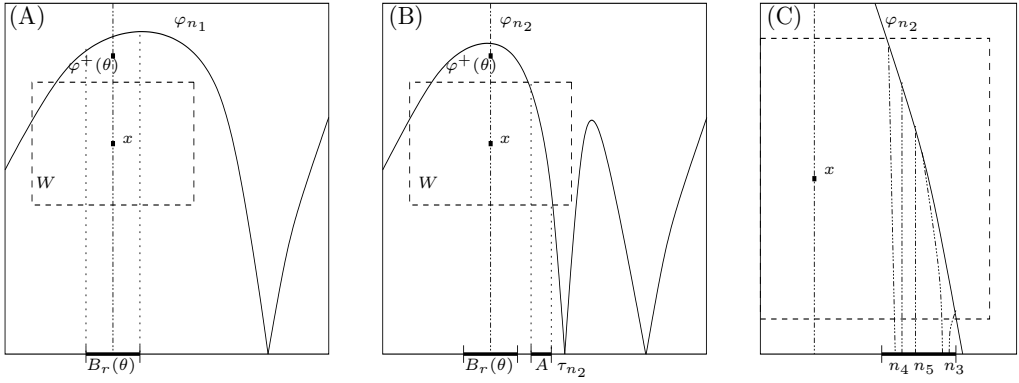, height=0.25\textheight,
    width=\textwidth} 
 \caption{\small 
   Taking the behaviour described above for granted we could argue as
   follows: (A)~Choose any $x \in [0,\varphi^+(\theta))$ and any small
   box $W$ around \thx\ which is below $\varphi^+(\theta)$. Fix
   $\neins \in \N$, such that the width of the peaks appearing at step
   $n_1$ or later is already small compared to $W$. As
   $\varphi_{\neins}(\theta) \geq \varphi^+(\theta)$ and due to the
   continuity of $\varphi_{\neins}$, this graph will be above $W$ on a
   whole neighbourhood $B_r(\theta)$ of $\theta$.  (B)~Now let $n_2$
   be the first time greater than $n_1$ a peak comes close to
   $B_r(\theta)$ (such that the value of the iterated upper boundary
   lines changes inside of $B_r(\theta)$). The end of that peak which
   is contained in $B_r(\theta)$ is still above $W$ (for now we ignore
   the case where the peak covers all of $B_r(\theta)$; this can be
   treated similarly). On the other
   hand $\varphi_{n_2}$ is pinched at $\tau_{n_2}$. Thus
   $\varphi_{n_2}$ crosses $W$ on an interval $A$, and as the slope of
   $\varphi_{n_2}$ is at most $\beta\alpha^{n_2}$ for some $\alpha >
   1, \ \beta>0$ this interval still has a certain length.  (C)~offers
   a closer look at what happens in the interval $A$ after $n_2$. As
   it takes a long time until the next peaks hit this interval and the
   width of these peaks decays exponentially, they will not cover all
   of $A$. (Only the first three of these peaks at times $n_3, n_4$
   and $n_5$ are depicted.) But this means $\varphi_{n_2}$ already
   coincides with $\varphi^+$ at least for some $\theta' \in A$.
   Therefore $\Phi^+ \cap W \neq \emptyset$, and as this works for any
   such box $W$ we can conclude $\thx \in \clPhiplus$.}
 \label{fig:heuristics}
\end{minipage}
\end{figure}

Figure \ref{fig:bounding} shows the first six iterates of the upper
boundary line. The pictures suggest a strikingly simple scenario:
Suppose $T$ is pinched exactly at $\theta = 0$ (i.e.\ $T_0 \equiv 0$
and all other fibre maps are striktly monotonically increasing) and
let $\tau_n := n\omega \bmod 1$. Then in the $n$-th iteration step a
new ``peak'' appears at $\tau_n$. The graphs $\varphi_{n-1}$ and
$\varphi_n$ differ on a small interval centered around $\tau_n$,
but apart from that they seem to coincide.  In addition, the peaks
seem to get steeper at an exponential rate, and accordingly the width
of the intervals decays exponentially. Figure \ref{fig:heuristics}
briefly explains on an heuristic level how such an observation can be
used to show that $\Phi^+$ is indeed dense in $[0,\varphi^+]$. A
slightly refined version of this will then determine our strategy in
the proof of Thm.\ \ref{thm:topcl} (see beginning of Section
\ref{Mainresult}).

Unfortunately, there are some complications in the general case.  When
the fibre maps are strictly monotonically increasing, the graphs
$\varphi_n$ will be strictly monotonically decreasing and thus two of
them can never exactly coincide (except on the orbit of the pinched
fibre). Still, it is possible to control the behaviour of these
graphs, but instead of concentrating on the peaks as our observations
suggest we must rather take care of the regions away from the peaks
(see Lemma \ref{lem:control}). However, this will turn out to be
sufficient for our purposes.

In order to obtain the required control we must put several
restrictions on the systems we study. As mentioned, for parameter
families as in (\ref{eq:pinchedskews}) these will basically translate
to a lower bound on the parameter. 

\ \\
Let $\tau_n := n\omega \bmod 1 \in \kreis$ as above and suppose
$\omega \in [0,1] \setminus \Q$ satisfies the diophantine condition
\begin{equation}
  \label{eq:diophantine}
  d(\tau_n,0) \ \geq \ c\cdot n^{-d}
\end{equation}
for some $c,d>0$. Further assume $T$ is a pinched skew product with
the following additional properties\footnote{For the sake of
  simplicity we state the conditions in terms of derivatives, although
  the system may not be differentiable everywhere. Then these
  conditions should always be understood in the way that the mentioned
  inequalities hold for the liminf and limsup of the respective
  difference quotients.}:
\begin{eqnarray}
  \label{eq:standingasmp}
  & \bullet & D\Tth(x) \leq \alpha \ \forall \thx  \textrm{ and }
  D\Tthx \leq \alpha^{-\gamma} \textrm{ whenever }  x
  \geq 1 \textrm{ for some } \alpha > 2, \ \gamma > 0. \hspace{3eM}
   \label{eq:exp/contr}\\
  & \bullet & |\frac{\partial}{\partial \theta} \Tthx | \leq \beta \
  \forall \thx \ \ 
  \textrm{for some} \ \beta > 0. \label{eq:beta} \\
  & \bullet & \textrm{Let } m\in\N \textrm{ satisfy } m \geq
  \label{eq:gamma} 
  4+\frac{4}{\gamma} \textrm{ and } \\
  & & \hspace{2eM} a \geq (m+1)^d \ \label{eq:a} , \\
  & & \hspace{2eM} b \leq  c \ \textrm{ with } \ d(\tau_n,0) \geq b \
  \forall n = 1 \ld m-1 \ . \label{eq:b} \\
  & & \textrm{Then suppose} \nonumber \\
  & & \hspace{1eM} \Tthx \geq \min\{1,ax\} \cdot \label{eq:reference}
  \min\{1,\frac{2}{b}d(\theta,0)\} \ \   
  \forall \thx \in \kreis \times I \hspace{1eM} \textrm{(reference
    system)  } 
\end{eqnarray}

\begin{bem} \label{bem:conditions}
  Here the reference system $R$ with fibre maps $R_\theta(x) =
  \min\{1,ax\} \cdot \min\{1,\frac{2}{b}d(\theta,0)\}$ is only used
  implicitly as a lower bound for the original system. However, it
  should be mentioned that this simplified system proved to be very
  useful in the development of the presented ideas. Actually all
  results were first derived in this simpler setting, where the
  heuristics described above can be directly converted into a rigorous
  proof.
  
  Note also that (\ref{eq:exp/contr}) and (\ref{eq:reference})
  together imply a certain ``expansion-contraction-property'': By
  (\ref{eq:reference}) the system is expanding close to the 0-line and
  by (\ref{eq:exp/contr}) it is contracting further away. As we will
  see, the expansion is responsible for the fact that the peaks become
  steeper whereas the contraction will be needed to control what
  happens away from the peaks.
\end{bem}
We start with two simple observations:
\begin{lem}
  \label{lem:returntimes}
\alphlist
\item $d(\tau_n,0) \leq b \cdot a^{-i} \ \follows \ n \geq
  a^{\frac{i}{d}}$
\item $|\varphi_n'| \leq \beta \cdot \alpha^n$.
\listend
\end{lem}
\proof 
\alphlist
\item $d(\tau_n,0) \leq b \cdot a^{-i} \
  \stackrel{(\ref{eq:diophantine})}{\follows} c \cdot n^{-d} \leq b
  \cdot a^{-i} \ \stackrel{(\ref{eq:b})}{\follows} n^{-d} \leq a^{-i}
  \ \follows \ n \geq a^{\frac{i}{d}}$
\item We prove $|\varphi_n'| \leq \beta \cdot (\alpha^n - 1)$ by
  induction. Note that 
  \begin{equation}
    \varphi_{n+1}'(\theta) =  DT_{\theta-\omega}(\varphi_n(\theta
    -\omega))\cdot \varphi_n'(\theta-\omega) +
    \frac{\partial}{\partial \theta} 
    T_{\theta-\omega}(\varphi_n(\theta-\omega)) \ .
  \end{equation}
For $n=1$ the statement is obvious as $\varphi_0' \equiv 0$, $\alpha > 
2$ and
$|\frac{\partial}{\partial \theta} \Tth(x)| \leq \beta$
((\ref{eq:beta})). If it is satisfied for some $n \geq 1$ we get (from 
(\ref{eq:exp/contr}) and (\ref{eq:beta}))
\[
    |\varphi_{n+1}'(\theta)| \leq \alpha \cdot
    |\varphi_n'(\theta-\omega)| + \beta \leq \alpha \cdot \beta \cdot
    (\alpha^n -1) + \beta \leq \beta \cdot (\alpha^{n+1} - 1) \ ,
\]
again using $\alpha \geq 2$.
\listend

\begin{figure}[h!]   
\noindent
\begin{minipage}[t]{\linewidth} 
   \epsfig{file=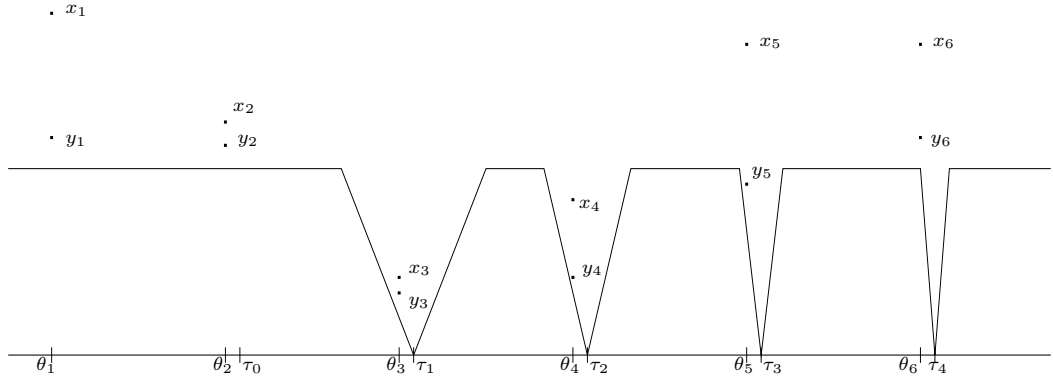, height=0.25\textheight,
    width=\textwidth} 
 \caption{\small Why do we have to expect
   $\varphi_{n-1}(\theta) - \varphi_n(\theta)$ to be rather large not
   only when $\theta \in J^{n-1}_n$, but also when it is in
   $J^{n-1}_j$ for some $j<n$. We illustrate this with a simple
   example where $n=6$ and $\theta=\theta_6$ is close to $\tau_4$.
   Ignoring the fact that we are on the circle, the first four peaks
   of the reference system are drawn in a straight order. They give a
   lower bound for the points $y_k := \varphi_{k-1}(\theta_k)$ and
   $x_k := \varphi_k(\theta_k) \ (k=1 \ld 6)$. In the first two steps 
   these points are above 1 and their distance is contracted But
   afterwards $\theta_k$ is inside of the $(k-2)$th peak, such that
   the points can be below 1 and their distance may be expanded again.
   Thus the assumption on $\theta$ in Lemma \ref{lem:control}(a) is
   needed to exclude that any of the $\theta_k$ is too close to
   $\tau_0=0$, maybe apart from the last $q-1$ steps.  The possibility
   that many different $\theta_k$ are close to $\tau_0$ is then ruled
   out by the diophantine condition on $\omega$.}
 \label{fig:control}

\end{minipage}
\end{figure}

\qed

\ \\
Let $J^n_j$ be an interval of length $b\cdot a^{-\frac{n}{m}}$
centered around $\tau_j$, i.e.\ 
 \begin{equation}
  \label{eq:intervals}
  J^n_j := (\tau_j - \frac{b}{2}\cdot a^{-\frac{n}{m}},\tau_j +
  \frac{b}{2}\cdot a^{-\frac{n}{m}}) \ .
\end{equation}
The following lemma now contains all the information about the
iterated upper boundary lines we need. The condition on $\theta$ is
probably be a little bit surprising at first: One might expect that
the difference between $\varphi_{n-1}$ and $\varphi_n$ must be very
small outside of a small interval around $\tau_n$. However, this is
not entirely true. As we shall see, it is only possible to ensure that
$|\varphi_{n-1}(\theta)- \varphi_n(\theta)|$ is small whenever
$\theta$ is sufficiently far away from ``most'' of the $\tau_j$, i.e.\ 
from $\tau_q,\tau_{q+1} \ld \tau_n$ at the same time, where $q$ is
relatively small compared to $n$. (The reason why we want to be able
to omit the first few $\tau_k$ is going to become obvious later, when
$q$ is specified in the proof of Thm.\ \ref{thm:topcl}.) The proof
below, together with Figure \ref{fig:control}, hopefully clarifies why
we have to admit $|\varphi_{n-1}(\theta) - \varphi_n(\theta)|$ to be
quite large not only when $\theta$ is close to $\tau_n$, but also when
it is close to $\tau_{n-1},\tau_{n-2},\ldots$ .
\begin{lem}
  \label{lem:control}
\alphlist
\item Let $\theta \notin \bigcup_{j=q}^n J^{n-1}_j$. If $q \leq
  \frac{n-1}{m}$, then 
  \begin{equation}
    \label{eq:control}
    |\varphi_{n}(\theta)-\varphi_{n-1}(\theta)| \ \leq L \cdot
    \alpha^{-(n-1)\lambda}
  \end{equation}
  where $\lambda :=  \gamma(1 - \frac{4}{m})-\frac{4}{m}$ is positive
  by (\ref{eq:gamma}). 
\item If $\theta \notin \bigcup_{j=q}^k J^{\tilde{n}-1}_j \cup
  \bigcup_{j = \tilde{n}+1}^\infty J^{j-1}_j$ and $q \leq
  \frac{\tilde{n}-1}{m}$ , then (\ref{eq:control}) equally holds for
  all $n \geq \tilde{n}$.
\item Let $\theta \notin \{ \tau_j \mid j \in \N \}$. Then there are
  infinitely many $n \in \N$ such that $\theta \notin \bigcup_{j=1}^n
  J^{n-2}_j$. 
\listend
\ \\
Note that the intervals used in (c) are slightly bigger than those
needed for the application of (a) to $\theta$. In fact, the difference 
is exactly $\frac{b}{2} \cdot (a^{\frac{n-2}{m\cdot d}} -
a^{\frac{n-1}{m\cdot d}})$ to either side. Therefore (a) can be
applied to all $\theta'$ from a small neighborhood of $\theta$.
\end{lem}
\proof 
\alphlist
\item Let $\theta_k:= \theta-n\omega+k\omega, \ x_{k} :=
  \varphi_k(\theta_k)$ and $s := \#\{1 \leq k < n \mid
  x_k < 1\}$. Now (\ref{eq:exp/contr}) implies
  \begin{eqnarray*}
    \lefteqn{\varphi_{n-1}(\theta) - \varphi_n(\theta) \ 
      = } \\
    & = &  (\varphi_0(\theta_1)-\varphi_1(\theta_ 1)) \cdot
    \prod_{k=1}^{n-1} 
     \frac{\varphi_k(\theta_{k+1})-
       \varphi_{k+1}(\theta_{k+1})}{\varphi_{k-1}(\theta_k) -
       \varphi_{k}(\theta_k)} \ \leq  \\ 
    & \leq & L \cdot \prod_{k=1}^{n-1}
     \underbrace{
       \frac{T_{\theta_k}(\varphi_{k-1}(\theta_k))  
       - T_{\theta_k}(\varphi_k(\theta_k))}
     {\varphi_{k-1}(\theta_k) - \varphi_k(\theta_k)}
      }_{\leq \ \alpha^{-\gamma} \ \textrm{whenever }
          \varphi_{k}(\theta_k) = x_k \geq 1 \ } 
     \ \leq \ L \cdot \alpha^{s - \gamma (n-1-s)} \
   \end{eqnarray*}
  Thus we are finished if we can show $s \leq \frac{4(n-1)}{m}$,
  because then $s - \gamma(n-1-s) \leq
  (n-1)(\frac{4}{m}-\gamma(1-\frac{4}{m})) = -(n-1)\lambda$. 
  
  In order to obtain an estimate on $s$, we first consider blocks of
  successive times where $x_k$ is smaller than 1. Inside of such a
  block $x_k$ is multiplied at least by the factor $a$ in each step,
  unless $\theta_k \in [-\frac{b}{2},\frac{b}{2}]$. Thus, most of the
  $x_k$ must even be below $\atel$. In fact, it is convenient to
  consider blocks $l+1 \ld p \ (0 \leq l < p < n)$ of $p-l$ successive
  times which satisfy $x_{l} \geq \atel$, \ $x_k < \atel \ \forall k =
  l+1 \ld p-1$, $x_{p} < 1$ and either $x_p \geq \atel$ or $x_{p+1}
  \geq 1$ or $p+1 = n$.  (In other words, either two blocks are
  seperated by at least one time where $x_k \geq 1$, or we start a new
  block when the threshold $\atel$ is reached.) Then
  necessarily $\theta_l \in [-\frac{b}{2},\frac{b}{2}]$ (otherwise
  $x_{l+1} \geq 1$), and the reference system (\ref{eq:reference})
  gives
  $x_{k+1} \geq a \cdot x_k \cdot \min\{1,\frac{2}{b} d(\theta_k,0)\}
  \ \forall k = l+1 \ld p-1$.  As $x_l \geq \atel$ and $x_p < 1$ we
  must have
  \[
     \prod_{k=l}^{p-1} \min\{1,\frac{2}{b}d(\theta_k,0)\} \ \leq
     a^{-(p-l-1)} \ .
  \]
  This means that $\sum_{1 \leq i < \infty} i \cdot \# \{ l \leq k < p
  \mid \frac{b}{2} \cdot a^{-i} \leq d(\theta_k,0) < \frac{b}{2}\cdot
  a^{-i+1} \}$ is an upper bound on the length $p-l$ of the block,
  such that every visit $d(\theta_k,0) \in [\frac{b}{2}\cdot
  a^{-i},\frac{b}{2}\cdot a^{-i+1})$ accounts for at most $i$ times
  $j>k$ where $x_j$ can be smaller that 1. Summing up the lengths of
  all blocks then gives an upper bound on $s$, namely
  \[
    s \leq \sum_{1 \leq i < \infty} i \cdot \# \{ 0 \leq k  < n-1 \mid
    \frac{b}{2} \cdot a^{-i} \leq d(\theta_k,0) < \frac{b}{2} \cdot
    a^{-i+1} \} \ .
  \]
  Now $\theta \notin \bigcup_{j=q}^n J^{n-1}_j$ implies
  \[
    d(\theta_k,0) \geq \frac{b}{2} \cdot
    a^{-\frac{n-1}{m}} \ \forall k = 0 \ld n-q ,
  \]
  anything else leads to the contradiction $\theta \in
  J^{n-1}_{n-k}$. In the worst case $x_{n-q+1} \ld x_{n-1}$ are all
  below 1, but this still allows the estimate
\begin{eqnarray*}
  \lefteqn{s \ \leq} \\
  & \leq &  \ q -1 + \sum_{1 \leq i \leq 
      \frac{n-1}{m}+1} i \cdot \# \{ 0 \leq k \leq n-q \mid
    \frac{b}{2} \cdot  
    a^{-i} \leq     d(\theta_k,0) < \frac{b}{2} \cdot a^{-i+1} \} 
     \ \leq  \\
  & \leq & \frac{n-1}{m} -1 + \sum_{1 \leq i \leq \frac{n-1}{m}+1} \# \{
    0 \leq k \leq  n-q \mid 
    d(\theta_k,0) < \frac{b}{2} \cdot a^{-i+1} \} \ \leq \\
  & \leq & \frac{2(n-1)}{m} + \sum_{2 \leq i \leq  \frac{n-1}{m}+1}
    1+\frac{n-1}{a^{\frac{i-1}{d}}} \ \leq \ \frac{3(n-1)}{m} + (n-1) \cdot
    \sum_{i=1}^{\infty} a^{-\frac{i}{d}} \
    \stackrel{\scriptstyle (\ref{eq:a})}{\leq} \ \frac{4(n-1)}{m}
\end{eqnarray*}
For the step from the third line to the last note that $\theta_k$ can
always visit the interval $(-\frac{b}{2} \cdot a^{-i+1},\frac{b}{2}
\cdot a^{-i+1})$ once, but a second visit at time $l$ implies
$d(\theta_k,\theta_l) = d(\theta_{l-k},0) < b\cdot a^{-i+1}$ and
thus $l-k \geq a^{\frac{i-1}{d}}$ by Lemma \ref{lem:returntimes}(a)
(or $l-k \geq m$ when $i = 1$ by (\ref{eq:b})).
\item is a direct consequence of (a).
\item Let $\nnull \in \N$ such that 
\begin{equation} \label{eq:nnull}a^{\frac{n-2}{m\cdot d}} \geq n+1
  \ \forall n \geq \nnull
\end{equation}
and suppose for some $k \geq \nnull$ we have $\theta \in
\bigcup_{j=1}^{k} J^{k-2}_j$.  Let $\tau_l$ be closest to $\theta$
of all points $\tau_1 \ld \tau_k$.  Then there exists a unique $n > k$
such that
$d(\theta,\tau_l) \in [b\cdot a^{-\frac{n-2}{m}},b\cdot
a^{-\frac{n-3}{m}})$. Now assume for some $j > l$ we have
$d(\tau_j,\theta) < \frac{b}{2} \cdot a^{-\frac{n-2}{m}}$. This
implies 
$d(\tau_j,\tau_l) = d(\tau_{j-l},0) < b \cdot a^{-\frac{n-2}{m}}$ and
therefore
$j \geq j - l \geq a^{\frac{n-2}{m\cdot d}} \geq n+1$ by Lemma 
\ref{lem:returntimes}(a) and (\ref{eq:nnull}). Thus $n$ satisfies 
our assumption.
\listend

\qed

\ \\


\section{The topological closure of the upper bounding graph}
\label{Mainresult} 

Let us first reformulate the result of the last section: Fix some
$\theta \in \kreis$ not contained in the forward orbit of the pinching
point. When $\theta$ is inside of the $k$-th peak, then the value of
$\varphi_n(\theta)$ may change significantly at time $n=k$ and for a
certain time afterwards, but then the movement of the iterated upper
boundary lines ``settles down'' in a small neighbourhood of $\theta$,
until the next peak visits that neighbourhood. Lemma
\ref{lem:control}(c) guarantees that there will always be such times
where the $\varphi_n$ have ``settled down'' close to $\theta$.

This can now be used to show that indeed $\clPhiplus = [0,\varphi^+]$,
using a slight modification of the strategy sketched in Figure
\ref{fig:heuristics}.  The peaks correspond to the intervals
$J^{n-1}_n$. $n_2$ will be chosen as before, i.e.\ the first time
after $n_1$ the interval $J^{n_2-1}_{n_2}$ hits the neighbourhood
$B_r(\theta)$. Now the main difference to Figure \ref{fig:heuristics}
is, that we cannot control $\varphi_{n-1}-\varphi_n$ on
$J_{n_2}^{n_2-1}$ until $n$ is bigger than $n_3 := m(n_2+1) + 1$ (see
Lemma \ref{lem:control}(a), where we will have to choose $q=n_2+1$).
Thus, instead of looking at the graph $\varphi_{n_2}$ and showing that
it coincides with $\varphi^+$ for some $\theta' \in B_r(\theta)$
(compare Fig.\ \ref{fig:heuristics}(B) and (C)), we have to look at
$\varphi_{n_3}$ and show that this graph is already very close to
$\varphi^+$ for some $\theta' \in B_r(\theta)$ (by controlling
$\sum_{n=n_3+1}^\infty (\varphi_{n-1} - \varphi_n)$ via Lemma
\ref{lem:control}(b)).

It may be worth mentioning that the conditions
(\ref{eq:exp/contr})--(\ref{eq:reference}) we put on our system only
had to be used to obtain statement (a) of Lemma \ref{lem:control} .
From now on, this lemma contains everything we need to know about our
system (apart from the pinched structure and the diophantine rotation
number, of course).

\begin{thm}
  \label{thm:topcl}
  Suppose $T$ satisfies all of the assumptions (\ref{eq:0line}),
  (\ref{eq:pinching}) and (\ref{eq:diophantine})--(\ref{eq:reference}).
  Then
\[
    \clPhiplus = [0,\varphi^+] \ .
\]
\end{thm}
\proof \\
Let $x \leq \varphi^+(\theta) - 3\epsilon$, $\delta>0$ and consider boxes
$W':= B_\delta(\theta) \times B_{2\epsilon}(x), \ W:= B_\delta(\theta)
\times B_\epsilon(x)$.
It suffices to show that $W' \cap \Phi^+ \neq \emptyset$ for
any such $\epsilon$ and $\delta$. 

First of all, we will fix some $\nnull \in \N$ which satisfies certain
conditions. Roughly speaking, these imply that for all $n \geq \nnull$
the intervals $J^{n-1}_j$ are small enough in comparison with $\delta$
and that the time $l$ we have to wait until $J^{n-1}_{j+l}$ hits
$J^{n-1}_j$ again is very large. All the conditions are obviously
satisfied when $\nnull$ is large enough. Apart from that the reader
should be advised not to wonder about these requirements too long in the
beginning, but to check each of them only at the time when it is
actually used below and the motivation for choosing it becomes
apparent.
\\
Let $\nnull \in \N$ satisfy
\begin{eqnarray}
  \label{eq:nnull1}
   \nnull & \geq  & m + 1 
   \\
   \label{eq:nnull2}
   \frac{b}{2} \cdot ( a^{-\frac{n-2}{m}} -   
   a^{-\frac{n-1}{m}})  & \leq & \frac{\delta}{2} \ \ \ \forall n \geq
   \nnull  
   \\
  \label{eq:nnull3}
    |J^{n-1}_j| \ = \ b \cdot a^{-\frac{n-1}{m}} & \leq & \frac{\delta}{2}
    \ \ \ \forall n \geq \nnull
    \\
  \label{eq:nnull4}
    a^{\frac{n-1}{m\cdot d}} & > & m(n+1) + 1 - n \hspace{2eM} \forall
    n \geq \nnull
    \\
   \label{eq:nnull5}
    \sum_{j \geq a^{\frac{n-1}{m \cdot d}}} b \cdot a^{-\frac{j-1}{m}} 
    & \leq &  \frac{\epsilon}{\beta} \cdot \alpha^{-m(n+1)-1}  \hspace{2eM}
    \forall n \geq \nnull
    \\
   \label{eq:nnull6}
    \sum_{j=n+1}^{\infty} L \cdot a^{-(j-1)\lambda} & \leq & \epsilon
    \hspace{2eM} \forall n \geq \nnull
\end{eqnarray}
For (\ref{eq:nnull5}) note that the left side decays
super-exponentially with $n$ whereas the right side only decays
exponentially.

\ \\
Now take some $\neins \geq \nnull$ which satisfies 
\begin{equation}
  \label{eq:neins}
  \theta
\notin \bigcup_{j=1}^{\neins} J^{\neins-2}_j \ .
\end{equation}
Such a $\neins$ always exists by Lemma \ref{lem:control}(c). As
pointed out that lemma, we can apply statement (a) of
it to all $\theta'\in B_r(\theta)$ where $r:=\frac{b}{2} \cdot (
a^{-\frac{\neins-2}{m}} - a^{-\frac{\neins-1}{m}}) > 0$, and this is
also possible for all $n \geq \neins$ as long as $B_r(\theta) \cap
J^{j-1}_j = \emptyset \ \forall j = \neins + 1 \ld n$. In addition we
can assume, by reducing $r$ further if necessary, that
$\varphi_{\neins}(\theta') > x+2\epsilon \ \forall \theta' \in
B_r(\theta)$ (note that $\varphi_{\neins}$ is continuous and
$\varphi_{\neins}(\theta) \geq \varphi^+(\theta) \geq x + 3\epsilon$).

Let $n_2$ be the first integer greater than $\neins$ such that
$J^{n_2-1}_{n_2}$ hits $B_r(\theta)$ and define $n_3 := m(n_2+1)+1$.
(\ref{eq:nnull2}) and (\ref{eq:nnull3}) ensure that $J^{n_2-1}_{n_2}
\subseteq B_\delta(\theta)$. Further, (\ref{eq:nnull4}) ensures
\begin{equation}
  \label{eq:ndrei}
  \theta' \notin \bigcup_{j=n_2+1}^{n_3} J^{j-1}_{j} \ \ \ \forall
  \theta' \in J^{n_2-1}_{n_2} \ ,
\end{equation}
because $J^{n_2-1}_{n_2} \cap J_j^{j-1} \neq \emptyset$ for some $j >
n_2$ implies $d(\theta_{j},\theta_{n_2}) \leq b \cdot
a^{\frac{-n_2+1}{m}}$, and by Lemma \ref{lem:returntimes}(a) and
(\ref{eq:nnull4}) this cannot happen if $j - n_2 \leq m(n_2 + 1) + 1 -
n_2 = n_3 - n_2$ .

Now first assume $\theta \in J^{n_2-1}_{n_2}$. As $\varphi_{n_3}$ is
pinched at $\theta_{n_2}$ and at the same time $\varphi_{n_3}(\theta)
\geq \varphi^+(\theta) \geq x+3\epsilon$ this graph will cross the box
$W$ from below to above when going from $\tau_{n_2}$ to $\theta$.
Thus, if we define
\begin{equation}
   A := \{ \theta' \in J^{n_2-1}_{n_2} \mid \varphi_{n_3}(\theta') \in
  B_\epsilon(x) \} \ , \nonumber
\end{equation}
we can use Lemma \ref{lem:returntimes}(b) to obtain that
\begin{equation}
  \label{eq:crossingW}
 |A| \ \geq \ \frac{2\epsilon}{\beta} \cdot \alpha^{-n_3} 
  \ .
\end{equation}
Let 
\begin{equation}
  \label{eq:lastingset}
  B := A \setminus \bigcup_{j=n_3+1}^{\infty} J_j^{j-1}
\end{equation}
$B$ is still a set of positive measure: Let $n_4$ be the first
time where $J_n^{n-1}$ hits $J_{n_2}^{n_2-1}$ again. Then $n_4 \geq
a^{\frac{n_2-1}{m\cdot d}}$ (again Lemma \ref{lem:returntimes}(a))
and therefore
\begin{equation*}
  | A \cap \bigcup_{j=n_3+1}^{\infty} J_j^{j-1}| 
  \ \leq \
  \sum_{j=n_4}^\infty |J_j^{j-1}| 
  \ \leq \ 
  \sum_{j \geq a^{\frac{n_2-1}{m \cdot d}}} b \cdot a^{-\frac{j-1}{m}} 
  \ \stackrel{(\ref{eq:nnull5})}{\leq} \
 \frac{\epsilon}{\beta} \cdot \alpha^{-n_3} 
\end{equation*}
where we used $n_3 = m(n_2+1)+1$ in the last step. Thus we have $|B| \geq
\frac{\epsilon}{\beta} \cdot \alpha^{-n_3} $. 

Combining (\ref{eq:ndrei}) and (\ref{eq:lastingset}) we see that it is
possible to apply Lemma \ref{lem:control}(b) with $q=n_2+1$ and
$\tilde{n} = n_3$ for every $\theta' \in B$. Thus the values of the
graphs $\varphi_n \ (n \geq n_3)$ do not differ to much anymore on
$B$, and by using $\varphi_n \searrow \varphi^+$ we get
\[
\varphi_{n_3}(\theta') - \varphi^+(\theta') \ = \ 
\sum_{n=n_3+1}^\infty \varphi_{n-1}(\theta') - \varphi_n(\theta') \ 
\leq \ \sum_{n=n_3+1}^\infty L \cdot a^{-(n-1)\lambda} \ 
\stackrel{(\ref{eq:nnull6})}{\leq} \ \epsilon \ \ \ \forall \theta'
\in B \ .
\]
As $\varphi_{n_3}(\theta') \in B_\epsilon(x) \ \forall \theta' \in B$
and $B \subseteq B_\delta(\theta)$ this proves $\clPhiplus \cap W'
\neq \emptyset$.

When $\theta$ is not contained in $J^{n_2-1}_{n_2}$ we can equally
conclude that $\varphi_{n_3}$ crosses $W$ inside of the interval
$J^{n_2-1}_{n_2}$, such that the above argument applies in exactly the
same way. To see this, note that when $\theta \notin J^{n_2-1}_{n_2}$,
then at least one endpoint of $J^{n_2-1}_{n_2}$ is contained in
$B_r(\theta)$. Denote it by $\bar{\theta}$. As $\bar{\theta} \notin
J_{n_2}^{n_2-1}$ (recall that these intervals are open) we can combine
(\ref{eq:neins}), the definition of $n_2$ and (\ref{eq:ndrei}) (which
also applies to $\bar{\theta}$ as the endpoint of $J^{n_2-1}_{n_2}$)
to obtain that $\bar{\theta} \notin \bigcup_{j=1}^{n} J^{n-1}_j \ 
\forall n = \neins \ld n_3$ and as $n_1 \geq m+1$ by (\ref{eq:nnull1}) 
we can apply Lemma \ref{lem:control}(a) with $q=1$ to obtain
\[
\varphi_{n_3}(\bar{\theta}) \ \geq \ \varphi_{n_1}(\bar{\theta}) - \!
\! \! \sum_{j=n_1+1}^{n_3} \varphi_{j-1}(\bar{\theta}) -
\varphi_j(\bar{\theta}) \ \geq \ x + 2\epsilon - \! \! \sum_{j =
  n_1+1}^{n_3} L \cdot a^{-(n-1)\lambda}  
\stackrel{(\ref{eq:nnull6})}{\geq}  x + \epsilon
\]
Again $\varphi_{n_3}(\theta_{n_2}) = 0$ and thus the graph has to
cross $W$ between $\theta_{n_2}$ and $\bar{\theta}$. This completes
the proof of the theorem.

\enlargethispage*{20pt}
\qed 

\ \\
We now want to apply this to parameter families $T=T_\alpha$
with fibre maps given by (\ref{eq:pinchedskews}). To that end suppose
$g: \kreis \ra \kreis$ and $f : \R^+ \ra \R^+$ are continuous
functions which satisfy the following assumptions:
\begin{eqnarray}
  \label{eq:gcond1}
  & \bullet & g \textrm{ is differentiable and strictly positive on }
  \kreis \setminus \{ \theta_0 \}. \hspace{10eM} \\
  \label{eq:gcond2}
  & \bullet & g(\theta_0) = 0 \hspace{2eM}  \\
  \label{eq:gcond3}
  & \bullet & \lim_{\theta \nearrow \theta_0} - g'(\theta) > 0 \
  \textrm{ and } \ \lim_{\theta \searrow \theta_0} g'(\theta) > 0 \\
  \label{eq:fcond1}
  & \bullet & f \textrm{ is differentiable and monotonically
    increasing.} \\
  \label{eq:fcond2}  
  & \bullet & f(0) = 0 \ , \ f'(0) > 0 \\
  \label{eq:fcond3}
  & \bullet &
  \alpha^{1+\gamma} f'(\alpha) \ra 0 \ (\alpha \ra \infty) \textrm{
    for some } \gamma>0  
\end{eqnarray}

\begin{cor}
  \label{cor:fgsystem} 
  Suppose $\omega$ satisfies a diophantine condition
  (\ref{eq:diophantine}) and $T$ is given by 
  \[
     \thx \mapsto (\thom,\alpha g(\theta)f(x))
  \]
  with maps $f$ and $g$ satisfying the above assumptions
  (\ref{eq:gcond1})--(\ref{eq:fcond3}). Then there exists an $\alpha_0
  > 0$ such that
\[
    \clPhiplus = [0,\varphi^+]
\]
whenever $\alpha \geq \alpha_0$.
\end{cor}
\proof \\
Choose $m \geq 4+\frac{4}{\gamma}$.  W.l.o.g.\ we can assume $\theta_0
= 0$. First note that the system
\begin{equation}
    \thx \mapsto (\thom,\alpha_1\alpha_2g(\theta)f(x))
\end{equation}
is topologically conjugated to the system
\begin{equation}
   \label{eq:conjsystem}
    \thx \mapsto (\thom,\alpha_1g(\theta)f(\alpha_2x))
\end{equation}
by $h : \thx \mapsto (\theta,\alpha_2x)$. (\ref{eq:gcond1}),
(\ref{eq:gcond3}) and $f'(0) > 0$ imply that this system is bounded
from below by a suitable reference system (\ref{eq:reference}) if
$\alpha_1$ and $\alpha_2$ are large enough. Once $\alpha_1$ is fixed,
(\ref{eq:fcond3}) implies that the expansion-contraction-condition
(\ref{eq:exp/contr}) is satisfied when $\alpha_2$ is sufficiently
large.  It is also easy to see that the system is pinched, has an
invariant $0$-line, and that $\frac{\partial}{\partial \theta} \Tthx$
is bounded for any $\alpha_1$ and $\alpha_2$. Thus all the assumptions
of Thm.\ \ref{thm:topcl} are satisfied, and the conjugacy $h$ of
course preserves the property we are interested in.

\qed

\begin{example}
  \label{ex:tanhfamily}
  Consider $g(\theta) = |\sin(\pi\theta)|$, $f(x) = \tanh(x)$ and let
  $\omega$ be the golden mean. The smallest possible value for $m$ in
  (\ref{eq:gamma}) is $m=5$ (if we can take $\gamma > 4$, which we
  verify below). The closest return up to time $m-1=4$ is $\omega^3
  \approx 0.236...$, which we take for $b$. Now we just try if $a=8$
  works. To that end note that $3|\sin(\pi\theta)| >
  \min\{1,\frac{2}{\omega^3}d(0,\theta)\}$ and $\frac{4}{3}\tanh(1) >
  1$. Thus $\Tthx = 4|\cos(\pi\theta)|\tanh(8x)$ satisfies
  (\ref{eq:reference}) with $a=8$ and $b=\omega^3$. As $\tanh'(x) =
  \cosh^{-2}(x) < 2e^{-2x}$, it is also easily checked that
  (\ref{eq:exp/contr}) holds for some $\gamma >4$ ($2e^{-16} <
  32^{-4}$). Thus we can choose $\alpha_0 = 32$ for this particular
  parameter family.
  
  Compared to this, the lower bound for the existence of an SNA
  obtained from (\ref{eq:snacondition}) is $\log \alpha >
  -\int_{\kreis} \log |\cos(\pi\theta)| d\theta$ or $\alpha >
  \exp\left(-\int_{\kreis} \log |\cos(\pi\theta)| d\theta\right) = 2$,
  still leaving a certain gap between the two conditions. Of course,
  by suitable modifications of the assumptions and in the proofs
  above, this gap could be closed a little bit further (especially in
  the case of the reference system mentioned in Rem.\ 
  \ref{bem:conditions}), but closing it completely that way seems to
  be out of reach. However, as the approach used here requires a
  rather strong control about the behaviour of the systems, it does
  not seem very surprising that ``something gets lost along the way''.
  Thus the existence of this gap does not necessarily have to mean
  that counterexamples should be expected for lower parameters
  $\alpha$.  After all, there does not seem to be a qualitative change
  in the behaviour of the iterated upper boundary lines for smaller
  $\alpha$ (consider Fig.\ \ref{fig:bounding} for example),
  as long as the upper bounding graph is not equal to 0 anyway.
\end{example}

\begin{bem}
  \label{bem:numerics}
  When looking at pictures from numerical simulations as in Figure
  \ref{fig:sna}, the result we obtained may seem a little bit
  surprising at first.  There, the larger we choose the parameter
  $\alpha$, the ``thinner'' the SNA seems to be around the 0-line.
  However, our approach offers a perfect explanation for this. If
  $\alpha$ is large, the peaks will become steeper and narrower very
  fast. Thus they will be too small to be detected numerically quite
  soon, and the SNA seems to coincide with the iterated upper boundary
  line after a few steps already. On the other hand, as we have seen
  it is exactly this fast decay of the width of the peaks which
  enabled us to prove our result.
\end{bem}

The following remark contains a rather informal discussion of two
possible slight generalizations of Thm.\ \ref{thm:topcl} and its
corollary. We refrained from including them in the statement of the
theorem, because although the basic idea does not change at all, this
would have made the proof far less readable. The most prominent
examples are already covered anyway.
\begin{bem} \label{bem:generalizations}
\alphlist
\item First, the diophantine condition on $\omega$ can be weakened to
  some extent. It was used to ensure that the return time to an
  interval of length $b \cdot a^{-i}$ grows exponentially in $i$ (see
  Lemma \ref{lem:returntimes}(a)), but actually all we need is that
  this quantity grows faster than linearly. This is the case as long
  as the rotation number $\omega$ satisfies $d(0,\tau_n) \geq s
  \cdot e^{-\sqrt[d]{n}}$ for some $s>0, d > 1$ (then the return times
  will asymptotically grow at least like like $i^d$). The requirements
  on the other parameters will certainly have to be stronger in this
  case, but apart from that the proof above needs only the slightest
  modifications.
\item Second, instead of working with the reference system
  (\ref{eq:reference}),
  it is also possible to use a system of the kind $R_\theta(x) =
  \min\{1,ax\} \cdot \min\{1,(\frac{2}{b}d(\theta,0))^p\}$, i.e.\ to
  replace the ``sharp peak'' at $\theta=0$ with a critical point of
  finite order $p$. Then, instead of counting the visits in the
  intervals $(-\frac{b}{2} a^{-i},\frac{b}{2}a^{-i})$ one would have
  to counts the visits in $(-\frac{b}{2}
  a^{-\frac{i}{p}},\frac{b}{2}a^{-\frac{i}{p}})$ (compare the proof of 
  Lemma \ref{lem:control}(a)). The diophantine condition (or the
  one mentioned in (a)) still guarantees that the resulting sums
  converge, only the conditions on the parameters must certainly be
  altered again. 
\listend
\end{bem}


\section{Upper bounding graphs which contain isolated points}
\label{Counterexamples} 

For the construction of counterexamples, we need some facts about the
combinatorics of irrational rotations, as they can be found in
\cite{demelo/vanstrien:1993} (for example). Let $\nfolge{\an}$ be the
coefficients of the continued fraction expansion of $\omega$,
\nfolge{\qn} the closest return times and \nfolge{\sigma_n} the
closest returns. Then the following equations hold
\begin{eqnarray}
    && \sigma_n  =  \qn \omega \bmod 1 = \tau_{\qn}  \label{eq:returns} \\
  && q_0 = 1 ,\ \ q_1 = a_1, \ \ q_{n+1}=a_{n+1}\qn + q_{n-1} 
  \label{eq:returntimes} \\
  && \sigma_0 = \omega, \ \ \sigma_1 = 1 - a_1\omega, \ \ \sigma_{n-1} =
  a_{n+1}\sigma_n + \sigma_{n+1}  \label{eq:returnsII} \\
  && \frac{1}{q_{n+2}} \leq \sigma_n \leq \frac{1}{q_{n+1}}
  \label{eq:returntimesII} 
\end{eqnarray}
For the remainder of this section we will assume that 
\begin{equation}
  \label{eq:growthcondition}
  \sum_{i=0}^\infty \frac{q_i}{q_{i+1}} \ < \ \infty \ .
\end{equation}
As $\frac{1}{2a_{n+1}} \leq \frac{q_i}{q_{i+1}} \leq
\frac{1}{a_{i+1}}$ (see (\ref{eq:returntimes})) this is equivalent to 
\begin{equation}
  \label{eq:ancondition}
  \sum_{i=1}^\infty \frac{1}{a_i} \ < \ \infty \ .
\end{equation}
\textbf{Construction of a suitable function $\mathbf{g}$:} Choose
$\neins \in \N$ such that $\sum_{i=\neins}^\infty \frac{q_i}{q_{i+1}} <
  \halb$ and, for the sake of simplicity, $\sigma_{\neins}$ is to the
  right of $0$ (in a local sense).  Then all points $\sigma_{\neins +
    2i} \ (i \in \N)$ will be to the right of 0, whereas all points
  $\sigma_{\neins + 2i + 1} \ (i \in \N)$ will be to the left. Let
\begin{eqnarray}
  \label{eq:Inintervals}
  I_n & := & \left\{ \begin{array}{ll} \ 
      [-\sigma_n+\sigma_{n+2},-\sigma_{n+2}] & \textrm{ if } n-\neins
      \textrm{ is even } \\ \ 
      [-\sigma_{n+2},-\sigma_n+\sigma_{n+2}]
 & \textrm{ if } n-\neins
      \textrm{ is odd } \end{array} \right. \\
  \label{eq:In'intervals}
  I_n' & := & \left\{ \begin{array}{ll} \ 
      [-\sigma_n,-\sigma_{n+2}] & \textrm{ if } n-\neins
      \textrm{ is even } \\ \ 
      [-\sigma_{n+2},-\sigma_n]
 & \textrm{ if } n-\neins
      \textrm{ is odd } \end{array} \right. 
\end{eqnarray}
and   
$I := \bigcup_{n=\neins}^\infty I_n' \cup \{ 0 \} =$
$[-\sigma_{\neins},-\sigma_{\neins+1}]$.
Now choose any $a > 2$ and a function $g$ with the property that
$g_{|I_n} \equiv a^{-\qn}$ and $g_{|I_n'\setminus I_n}$ joins the two
levels $a^{-\qn}$ on $I_n$ and $a^{-q_{n-2}}$ on $I_{n-2}$ in a 
continuous and monotone way (such that 
$g(I_n') = [a^{-\qn},a^{-q_{n-2}}])$.

\ \\
\begin{figure}[h!]   
\noindent
\begin{minipage}[t]{\linewidth} 
   \epsfig{file=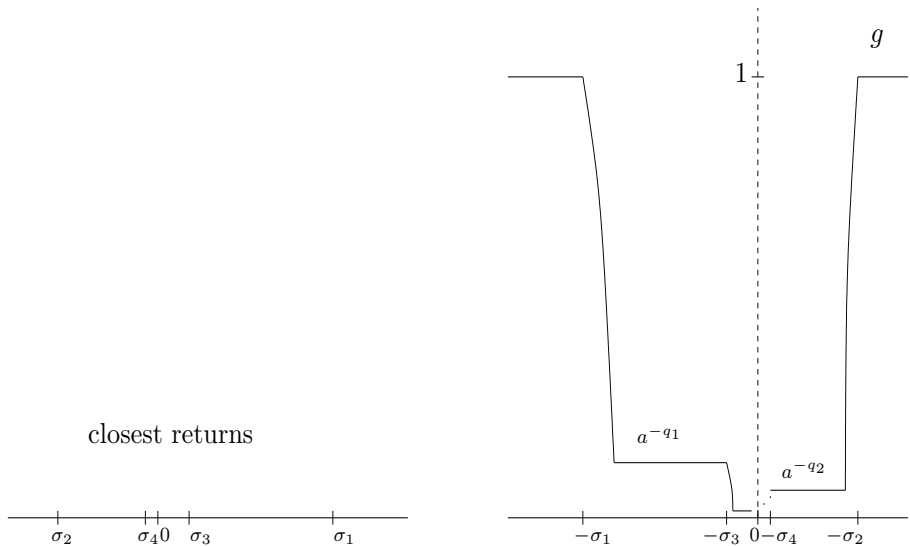, height=0.35\textheight,
    width=0.9\textwidth} 
 \caption{\small Construction of the function $g$ (with $\neins = 1$).}
 \label{fig:gconstruction}
\end{minipage}
\end{figure}
\ \\ 
Note that $|\log g|$ is integrable:      
\[ 
\int_{\kreis} |\log g(\theta)| \ d\theta \leq \log a \cdot
\sum_{i=\neins}^\infty q_i \cdot |\sigma_i|
\stackrel{(\ref{eq:returntimesII})}{\leq} \log a \cdot
\sum_{i=\neins}^\infty \frac{q_i}{q_{i+1}} \leq \halb \log a \ .
\]
Let further $f(x) := \min\{1,ax\}$ and consider the system $T$ given 
by
\begin{equation}
  \label{eq:examplesystem}
  T : \thx \mapsto (\thom,g(\theta)f(x)) \ .
\end{equation}
\textbf{Claim 1:} \[ \varphi^+(0) \ = \ 1 \]
\proof \\
Note that $\varphi_n(\theta) \neq \varphi_{n-1}(\theta)$ requires
$\varphi_{n-1}(\theta-\omega) \neq \varphi_{n-2}(\theta-\omega)$, and
in our particular example also $\varphi_{n-1}(\theta-\omega) < \atel$.
Let $n$ be the first time that $\varphi_n(0) < 1$. Then necessarily
$\varphi_{n-j}(\tau_{-j}) < \atel \ \forall j=1 \ld n-1$. Thus for the
first $n-1$ steps the
backwards orbit of $(0,\varphi_n(0))$ is always in the expanding region 
and therefore
\[
  \varphi_{n}(0) = T_{\tau_{-n}}^n(1) = a^{n-1} \cdot
  \prod_{j=1}^n g(\tau_{-j}) \ .
\]
As $\varphi_n(0) < 1$ this implies
\[
\prod_{j=1}^n g(\tau_{-j}) \leq a^{-(n-1)} \ .
\]
However, we can estimate $-\log_a \prod_{j=1}^n g(\tau_{-j})$ simply
by counting how often $(\tau_{-j})_{j=1\ld n}$ visits the intervals
$I_n'$. As the time between two such visits must be at least $q_{n+1}$ 
($|I_n'| < |\sigma_n|$), this gives
\[
   -\log_a \prod_{j=1}^n g(\tau_{-j}) \leq \sum_{i=\neins}^\infty q_i
   \cdot \frac{n}{q_{n+1}} \leq \frac{n}{2} \ .
\]
W.l.o.g.\ we can assume that $n \geq 2$ ($\tau_{-1} \notin I$) and
thus arrive at a contradiction. This proves claim 1.

\ \\
\textbf{Claim 2:}
\[ \varphi^+(\theta) \leq \atel \ \forall \theta \in I \setminus \{ 0
\} \]
\proof \\
This follows directly from the properties of $g$: When $n-\neins$ is
even we have
\[ 
   \varphi_{\qn} \leq \atel \ \textrm{ on } \
   [-\sigma_n+\sigma_{n+2},0] + \qn\omega = [\sigma_{n+2},\sigma_n]
\]
and
$\bigcup_{i=0}^\infty [\sigma_{\neins+2i+2},\sigma_{\neins+2i}] =$
$(0,\sigma_{\neins}]$.
The same argument applies to the left side, which proves claim 2.
\ \\
\ \\
The two claims together imply that $(0,\varphi^+(0)) = (0,1)$ is
isolated in $\overline{\Phi^+}$, and by continuity of $T$ the same is
true for all points in the backwards orbit of this point.

We close with some final remarks concerning the regularity of the  
examples just constructed. First of all, it is possible to replace the
rather degenerate function $f$ by a differentiable and strictly 
monotonically increasing function $\tilde{f}$. To that end, suppose
$\tilde{f} : [0,2] \ra [0,2]$ satisfies
$\tilde{f}_{|[0,\atel]} \equiv f_{|[0,\atel]}$, $\tilde{f}' > 0$ and
$\tilde{f}(2) \leq 2$. Then the system 
\[
    \tilde{T}: \thx \mapsto (\thom,g(\theta)\tilde{f}(x))
\]
is bounded from below by (\ref{eq:examplesystem}), which ensures
$\tilde{\varphi}^+(0) \geq 1$. On the other hand, the same arguments
as in the proof of claim 2 show that
$\tilde{\varphi}^+(\theta) \leq \frac{2}{a} < 1 \ \forall$
$\theta \in [\sigma_{\neins+1},\sigma_{\neins}]$.
As $\log g$ is integrable (see above) $\tilde{T}$ also
fits perfectly into the framework of Section \ref{Pinchedskews}. 

Finally, although (\ref{eq:ancondition}) certainly rules out the
golden mean or any rotation number with bounded coefficients $\an$, it
still allows $\omega$ to be of diophantine type. Then it is possible to
show, by combining the diophantine condition (\ref{eq:diophantine})
with some elementary estimates obtained from
(\ref{eq:returntimes})--(\ref{eq:returntimesII}), that
$\frac{a^{-q_{n-2}}}{\theta_{n+2}}$ converges to zero as $n$ goes to
infinity. Thus the connecting pieces of $g$ on $I_n' \setminus I_n$
can be chosen such that $g'(\theta) \ra 0$ as $\theta \ra 0$, in which
case $g$ is even differentiable. However, more or less the same
estimates yield that $g$ will then be infinitely flat at $\theta=0$,
such that it cannot be bounded from below by any suitable reference
system (compare Rem.\ \ref{bem:generalizations}(b)).

%
%
%

\bibliography{snaphysics,qpfs} \bibliographystyle{alpha}

\end{document}